\documentclass[12pt]{article}
\title{On $\alpha$-largeness and the Paris--Harrington principle in $\RCA$ and $\RCAst$}
\author{Florian Pelupessy \\ Mathematical Institute,
Tohoku University}
\date{}

\usepackage[english]{babel}
\usepackage{amsmath, amssymb, mathtools}
\usepackage{url}
\usepackage{xcolor}
\usepackage{hyperref}
\usepackage{tikz}
\usepackage{multicol}
\usepackage{parskip}
\usepackage{makeidx}
  \makeindex

 \usepackage[initials]{amsrefs}

\BibSpec{article}{%
  +{}{\PrintAuthors}				{author} %
  +{,}{ \textit}					{title}
  +{,}{ }							{journal}
  +{}{ \textbf}					{volume}
  +{,}{ }							{note}
  +{}{ \parenthesize}				{date}
  +{:}{ }							{pages}
  +{.}{}							{transition}
  +{}{ }							{review}
  +{\\}{ \url }							{url}
}

\BibSpec{thesis}{%
  +{}{\PrintAuthors}				{author} %
  +{,}{ \textit}					{title}
  +{,}{ }							{journal}
  +{}{ \textbf}					{volume}
  +{,}{ }							{note}
  +{,}{ }							{date}
  +{:}{ }							{pages}
  +{.}{}							{transition}
  +{}{ }							{review}
  +{\\}{ \url }							{url}
}

\usepackage{fancyhdr}
\makeatletter
\usepackage[T1]{fontenc}
\usepackage{libertine}

\newtheorem{theorem}{Theorem}
\newtheorem{lemma}[theorem]{Lemma}
\newtheorem{definition}[theorem]{Definition}

\newtheorem{corollary}[theorem]{Corollary}

\newenvironment{remark}[1][Remark]{\begin{trivlist}
\item[\hskip \labelsep {\bfseries #1}]}{\end{trivlist}}

\newcommand{\qed}{
\begin{flushright}
$\Box$
\end{flushright}
}

\newcommand{\conc}{%
  \mathord{
    \mathchoice
    {\raisebox{1ex}{\scalebox{.7}{$\frown$}}}
    {\raisebox{1ex}{\scalebox{.7}{$\frown$}}}
    {\raisebox{.7ex}{\scalebox{.5}{$\frown$}}}
    {\raisebox{.7ex}{\scalebox{.5}{$\frown$}}}
  }
}

\newcommand{\RCA}{\mathrm{RCA}_0}
\newcommand{\EFA}{\mathrm{EFA}}

\newcommand{\N}{\mathbb{N}}

\newcommand{\RT}{\mathrm{RT}}

\newcommand{\CNF}{\mathrm{CNF}}
\newcommand{\lh}{\mathrm{lh}}
\newcommand{\MC}{\mathrm{MC}}

\newcommand{\RCAst}{\mathrm{RCA}_0^{\displaystyle{*}}}

\begin{document}
\maketitle
\begin{abstract}%
We examine, within $\RCA$, the treatment by Ketonen and Solovay on the use of $\alpha$-largeness for giving an upper bound for the Paris--Harrington principle. This proof works fine in $\RCAst$ for every fixed standard dimension. We also show how to modify the arguments to work within $\RCAst$ for unrestricted dimensions. To the author's knowledge, this is the first time that it is confirmed that the treatment can be done within $\EFA$ without some transfinite induction added.
\end{abstract}


\section{Introduction}
In \cite{ketonensolovay} Ketonen and Solovay prove the following theorem:

\begin{theorem}[Ketonen--Solovay]\label{thm:ketonensolovay}
If $X\geq  3$ is $\omega_{d+1}(c+5)$-large, then every colouring $C\colon [X]^{d+2} \rightarrow c$ has homogeneous $H\subseteq X$ of size $> \min H$.  
\end{theorem}
This note consists of a modified version of the original presentation which should be better suited for reading from the reverse mathematics viewpoint. This may be of interest in light of the theorem's use, by Patey and Yokoyama, in the conservativity result for $\RT^2_2$ in  \cite{pateyyokoyama}. As stated there, once it is understood, the original proof, for $d=0$ and standard $c$ (Lemma~\ref{lemma:ketonensolovayd=0} in this note), is not hard to be seen to be formalisable in $\RCA$, since one can restrict the uses of transfinite induction to transfinite induction on $\omega^{c+4}$. However, for readers unfamiliar with the subject matter, it is somewhat tedious to check this due to the distribution of the proof throughout the paper. 

Thanks to arithmetic conservativity (Corollary~IX.1.11 from \cite{sosoa}), this is also a confirmation that the Ketonen--Solovay theorem is provable within $\mathrm{I}\Sigma_1$, as asked for in Problem~3.37 in \cite{hajekpudlak}. In the case of fixed standard $d$, one can easily weaken the base theory to $\RCAst$ without modifying the proof. This shows that Ketonen and Solovay's copious use of transfinite induction on ordinals is readily circumvented for the theorem in question. 

Finally, we confirm that one can also weaken the base theory to $\RCAst$ for unrestricted $d$. Thanks to $\Pi^0_2$-conservativity (Corollary~4.9 in \cite{simpsonsmith}), this implies that the Ketonen--Solovay theorem is provable in elementary function arithmetic, $\mathrm{EFA}$.

The presentation within $\RCA$ is suitable for advanced master level students and those who are unfamiliar with the Ketonen--Solovay paper \cite{ketonensolovay}. We assume only basic knowledge on reverse mathematics in $\RCA$ as in II.1-II.3 from \cite{sosoa}. At some places we favour an intuitive description and we leave many of the details as exercises for the reader. The changes compared to Ketonen--Solovay are concentrated in Section~\ref{section:411replacement}, with the rest of the proof, in Section~\ref{section:proof}, being only slightly modified from the originals. If one is only interested in the case of $d=0$ (dimension 2, as used in \cite{pateyyokoyama}), the presentation ends at the remark  after Lemma~\ref{lemma:ketonensolovayd=0}.

In the last section we will describe how to modify our arguments to work within $\RCAst$.

\section{Ordinals below $\varepsilon_0$ in $\RCA$}\label{section:ordinals}
We will define the ordinals below $\varepsilon_0$ within $\RCA$ as in Definition 2.3 in \cite{simpson1988}.
\begin{definition}\label{def:ordinal}
We define the set $\mathcal{E}$ of notations of ordinals $<\varepsilon_0$ and order $<$ on $\mathcal{E}$ as follows:
\begin{enumerate}
  \item If $\alpha_0 \geq \dots \geq \alpha_n \in \mathcal{E}$, then $\omega^{\alpha_0} + \dots + \omega^{\alpha_n} \in \mathcal{E}$.
  \item $\omega^{\alpha_0} + \dots + \omega^{\alpha_n} < \omega^{\beta_0} + \dots + \omega^{\beta_m}$ if and only if:
  \begin{enumerate}
    \item $n<m$ and  $\alpha_i=\beta_i$ for all $i \leq n$, or:
    \item there is $i\leq \min\{ n,m\}$ with $\alpha_j=\beta_j$ for all $j < i$ and $\alpha_i < \beta_i$.
  \end{enumerate}
\end{enumerate}
We use $0$ to denote the empty sum, $0 < \alpha$ for all $\alpha \neq 0$, $1=\omega^0$, $n=\overbrace{1 + \dots +1}^{n}$, $\omega=\omega^1$, $\omega_0(\alpha)=\alpha$, $\omega_{d+1}(\alpha)=\omega^{\omega_d(\alpha)}$ and $\omega_d=\omega_d(1)$. 
\end{definition}
As usual, if $\alpha_n=0$, then $\alpha$ is called a successor, otherwise, when not equal to $0$, it is called a limit. One can define primitive recursive functions for ordinal-addition, natural (Hessenberg) sum, and ordinal multiplication on $\mathcal{E}$. Recall that, for $\alpha$ and $\beta$ as in Definition~\ref{def:ordinal}, the natural sum is:
\[
\alpha \oplus \beta = \omega^{\gamma_0} + \dots + \omega^{\gamma_{m+n+1}},
\]
where the $\gamma_i$'s are all the $\alpha_i$'s and $\beta_i$'s in descending order. The natural sum has the important property that none of the terms are lost, which can happen with ordinal addition. For example:
\[
\omega+\omega^2=\omega^2\neq \omega^2+\omega=\omega \oplus \omega^2.
\]
Every ordinal in $\mathcal{E}$ has a Cantor Normal Form:
\[
\alpha=_{\mathrm{CNF}}\omega^{\alpha_0} \cdot a_0 + \dots + \omega^{\alpha_n} \cdot a_n,
\]
where the $a_i$'s are positive integers and $\alpha_0 > \dots > \alpha_n$.
\begin{definition}[Maximal coefficient]
$\mathrm{MC}(0)=0$ and, given $\alpha=_{\CNF}\omega^{\alpha_0} \cdot a_0 + \dots + \omega^{\alpha_n} \cdot a_n>0$:
\[
\mathrm{MC}(\alpha)=\max \{ a_i , \mathrm{MC}(\alpha_i) \}.
\]
\end{definition}
\begin{definition}[Fundamental sequence]
For $\alpha = \omega^{\alpha_0} + \dots + \omega^{\alpha_n} \in \mathcal{E}$ and $x \in \N$, take $0[x]=0$, $(\alpha+1)[x]=\alpha$, and:
\begin{enumerate}
  \item If $\alpha_n=\beta+1$, then $\alpha[x] = \omega^{\alpha_0} + \dots + \omega^{\alpha_{n-1}} + \omega^\beta \cdot x$,
  \item If $\alpha_n$ is a limit, then $\alpha[x] = \omega^{\alpha_0} + \dots + \omega^{\alpha_{n}[x]}$.
\end{enumerate}
\end{definition}
\begin{definition}
A finite set $X=\{x_0 < \dots < x_{|X|-1}\}$ is called $\alpha$-large if:
\[
\alpha[x_0] \dots [x_{|X|-1}]=0.
\]
Any $X$ is $0$-large. 
\end{definition}
Any $\omega$-large set $X$ has size $> \min X$. Unless otherwise specified, we will assume $\alpha$-large sets to be strictly above $2$. 

At first glance one may think that we require transfinite induction to demonstrate properties of the fundamental sequences and $\alpha$-large sets. In the remainder of this section we avoid this usage to treat some properties for later use. 
\begin{lemma}\label{lemma:fundbound}
If $\omega_d >\alpha > \beta$ and $x>\MC(\beta)$, then $\alpha[x] \geq \beta$, where the inequality is strict if $\alpha$ is a limit. 
\end{lemma}
\emph{Proof:} Induction on $d$.
\qed
\begin{lemma}\label{lemma:shiftlarge}
For any $\alpha$ and any $x_0 < \dots < x_R$ , $\MC(\alpha)<y_0 < \dots < y_R$ with $0<x_i \leq y_i$ for all $i\leq R$: if  $\alpha[x_0] \dots [x_R]>0$, then $\alpha[y_0] \dots [y_R]>0$. 
\end{lemma}
\emph{Proof:} Use induction on $R$ with the aid of Lemma~\ref{lemma:fundbound} to show $\alpha[y_0] \dots [y_R] \geq \alpha[x_0] \dots [x_R]$.
\qed
\begin{lemma}\label{lemma:largerordinal}
For any $\alpha>\beta>0$ and any $\MC(\beta) < x_0 < \dots < x_R$,  we have that $\alpha[x_0] \dots [x_R]> \beta[x_0] \dots [x_R]$.
\qed
\end{lemma}
\emph{Proof:} Use induction on $R$ with the aid of Lemma~\ref{lemma:fundbound}.

Define the thrice iterated exponential: 
\[
E(x)=2^{2^{2^x}}.
\]

One can check that:
\begin{enumerate}
  \item The smallest $\omega$-large interval which contains $x$ as its minimal element is $[x, 2x]$.
  \item For $\omega^2$ this is bigger than $[x, 2^x\cdot x]$. 
  \item If $x \geq 3$, then $\omega^3[x] \dots [E(x)+x+8] > 0$.
\end{enumerate}  
  
The following lemma shows that $\omega^3$-large sets $X$ are larger than $E(\min X)$. This is a rather weak lower bound, since Ketonen and Solovay showed in their original proof that one can use the tower function instead of $E$. 
\begin{lemma}\label{lemma:Ebounds}
For any $3 \leq x_0 < x_1 < \dots$ we have $\omega^3[x_0] \dots [x_{E(x_0) + 8}] >0$.
\end{lemma}
\emph{Proof:} This follows from item (3) directly above and Lemma~\ref{lemma:shiftlarge}. 
\qed

\section{Theorem 4.11 replacement:}\label{section:411replacement}

Take:
\[
\Phi(\alpha)= \omega^3 \cdot \alpha + \omega^3 +l +2.
\]

\begin{lemma}[Theorem 4.11-Replacement]\label{lemma:411replacement}
Suppose that $2<X=\{x_0, \dots x_{|X|-1}\}$ is $\Phi(\gamma_0)$-large and $\gamma_0 > \gamma_1 > \dots > \gamma_j$  is such that $\mathrm{MC}(\gamma_{i}) \leq E(x_i+l)$. Then $j \leq |X|-1$.
\end{lemma}
The outline of the proof is as follows: Take $\alpha_0=\Phi(\gamma_0)$ and $\alpha_{i+1}=\alpha_i [x_i]$. By  $\Phi(\gamma_0)$-largeness we know that $\alpha_{|X|}=0$. In the original proof of Theorem 4.11 in \cite{ketonensolovay} it is shown that the $\Phi(\gamma_i)$'s are a subsequence of the $\alpha_i$'s. We will show that the  $\alpha_i$'s contain a subsequence whose $i$th elements are \emph{larger} than the corresponding $\Phi(\gamma_i)$'s:
\begin{center}
\begin{tabular}{ccccccccc}
$\alpha_{a_0}$ 			& $>$ 	& $\dots$ 	&$>$ 	& $\alpha_{a_i}$	& $>$ 	& $\dots$ 	& $>$ & $0$ \\
					&		&			&		& $\vee$			& & & & \\
$\Phi(\gamma_0)$ 	& $>$	& $\dots$ 	&$>$	& $\Phi(\gamma_i)$ 		& $>$	 	& $\dots$ & & \ \ ,
\end{tabular}
\end{center}

thus demonstrating the conclusion of the lemma. The core of this lemma, namely pointing out the subsequence which has this property, is contained in the \emph{Claim}. 

\noindent\emph{Proof:} Notice that: 
\begin{enumerate}
  \item $\mathrm{MC}(\omega^3 \cdot \alpha+x) \leq \max \{\mathrm{MC}(\alpha)+3, x \}$,
  \item $\mathrm{MC}(\alpha[x]) \leq \max( \mathrm{MC}(\alpha) , x )$,
  \item $E(x+1) > E(x) + 4$,
  \item Take $\beta_0 = \omega^3 \cdot \beta+ \omega^3$ and $\beta_{k+1}=\beta_k[x_{i+k}]$, then, by Lemma~\ref{lemma:Ebounds}: $\beta_{E(x_i)+8} > \omega^3 \cdot \beta$.
\end{enumerate}

Take:
\[
a_i=\left\{ 
\begin{array}{ll}
0  & \textrm{ if $i=0$}, \\ 
E(x_{i+l+1}) & \textrm{otherwise}.
\end{array}
\right.
\]

\noindent\emph{Claim:} $\alpha_{a_i} > \Phi (\gamma_i)$ for all $0<i \leq j$. 

\noindent\emph{Proof of the claim:} Induction on $i$. We show both the case $i=1$ and the induction step simultaneously. 

\noindent We have the following , if $i=0$ by notice~(4), otherwise by induction hypothesis and all four notices: 
\[
\alpha_{a_i + E(x_{a_i+l+1})+l+8} > \omega^3 \gamma_i,
\]
Therefore, thanks to $\gamma_{i+1} < \gamma_i$:
\[
\alpha_{a_i + E(x_{a_i+l+1})+l+8} > \omega^3 \gamma_i \geq \omega^3 \cdot (\gamma_{i+1} + 1)=\omega^3 \cdot \gamma_{i+1}+ \omega^3.
\]
Hence:
\[
\alpha_{a_{i+1}}  \geq \alpha_{a_i + E(x_{a_i+l+1})} > \Phi(\gamma_{i+1}),
\]
thus ending the proof of the claim, hence the lemma.
\qed
\begin{remark}
As a side note, the claim in the lemma also implies that the $\Phi(\gamma_i)$'s are a subsequence of the $\alpha_i$'s by using the following fact which can be shown using induction on $i$: 

If $\alpha_{j-i-1} > \beta \geq \alpha_j$ and $x_{j-i-1} > \MC(\beta)$, then $\beta=\alpha_l$ for some $j-i\leq l \leq j$.
\end{remark}






\section{How to prove the Ketonen--Solovay theorem}\label{section:proof}

We proceed with, essentially, the proofs from  Section 5 and 6 of \cite{ketonensolovay}. The proofs have been streamlined into our setting, defining trees as sets of sequences, as is usual in reverse mathematics. The use of tree arguments for proving Ramsey-type theorems is attributed to Erd\"os and Rado. The outline is as follows: 

\begin{enumerate}
\item
Show that, if $X$ is $(\omega\cdot c)$-large, then for every colouring $C\colon X \rightarrow c$ there exists $\omega$-large $C$-homogeneous $H \subseteq X$. This is Lemma~\ref{lemma:php}.
\item
Given $C\colon [X]^{d+2}\rightarrow c$, construct Erd\"os--Rado trees $T_i$ from the first $i$ elements of $X$. Derive, from these trees, a decreasing sequence of ordinals of length $|X|$. Use Lemma~\ref{lemma:411replacement} to determine that, if $X$ is ``large enough'' compared to $\alpha$, then $T_{|X|}$ contains an $\alpha$-large branch $Y$ such that the value of $C(x)$ on $[Y]^{d+2}$ does not depend on $\max x$. The case $d=0$ is handled in Lemma~\ref{lemma:ketonensolovayd=0}, the case $d>0$ is treated in Lemma~\ref{lemma:stepup}. 
\item
Using induction on $d$, derive Theorem~\ref{thm:ketonensolovay} from the above. 
\end{enumerate}

\begin{lemma}\label{lemma:php} If $X$ is $(\omega \cdot c)$-large then for every colouring $C\colon X \rightarrow c$ there exists $\omega$-large $C$-homogeneous $H \subseteq X$.
\end{lemma}
\emph{Proof:}
Since $X$ is $(\omega \cdot c)$-large it is the disjoint union of $\omega$-large sets: $X=X_0 \cup \dots \cup X_{c-1}$. Assume, without loss of generality, that the $\min C^{-1}(i)$'s are increasing. Assume, for a contradiction, that no $C^{-1} (i)$ is $\omega$-large. By induction on $i < c-1$ we have:
\[
\min C^{-1} (i) \leq \min X_i \ \& \ |\bigcup_{j \leq i} C^{-1} (j)| < |\bigcup_{j\leq i}X_j|,
\]
the latter being implied by the first, as the $X_i$'s are $\omega$-large whilst the $C^{-1}(i)$'s are not.
This implies  $|\bigcup_{j \leq c-1} C^{-1} (j)| < |X|$, a contradiction.
\qed

\begin{definition}
For $0<i\leq d+1$, $C\colon [X]^{d+1} \rightarrow c$, we say $Y\subseteq X$ is ${\min}_i$-$C$-homogeneous if the value of $C$, on $[Y]^{d+1}$, depends only on the first $i$ elements of its input:
\[
C(x_0, \dots , x_{i-1} , y_i , \dots y_d)= C(x_0, \dots , x_{i-1} , z_i , \dots , z_d)
\]
for all $x_0 < \dots < x_{i-1} < y_i < \dots < y_d$, $x_{i-1} < z_i < \dots < z_d$ from $Y$.

\end{definition}

\begin{lemma}\label{lemma:ketonensolovayd=0}
If $X$ is $(\omega^{c+3}+\omega^3 + c+4)$-large, then every colouring $C\colon [X]^2 \rightarrow c$ has homogeneous $H\subseteq X$ of size $> \min H$.  
\end{lemma}
\emph{Proof:} Given $X=\{ x_0 < \dots < x_{|X|-1}\}$ and $C\colon [X]^{2} \rightarrow c$, by the previous lemma it is sufficient so show that $X$ has an $(\omega\cdot c)$-large  ${\min}_1$-$C$-homogeneous subset. Assume, for a contradiction, that is not the case.

Define $T_0 \subset \dots \subset T_{|X|}$ as follows: $T_0= \{ \emptyset \}$ and 
\[
T_{i+1}= T_i \cup \{ \sigma \conc x_i\},
\]
where $\sigma \in T_i$ is the leftmost of maximum length such that $\{ \sigma_0 ,  \dots, \sigma_{\mathrm{lh}(\sigma )-1} ,x_i\}$ is ${\min}_1$-$C$-homogeneous. By construction, if $\sigma\conc y , \sigma\conc z \in T_i$ then:
\[
C( \sigma_{\lh(\sigma)-1} , y ) \neq C(\sigma_{\lh(\sigma)-1}, z).
\]
So the number of branches of $\sigma$ has upper bound $c$. 

Let $(\omega\cdot c) [\sigma_0] \dots [\sigma_{\lh(\sigma)-1}]=\omega \cdot d_\sigma + r_\sigma >0$

Define: $n_{\sigma, i} = (c+1)^{r_\sigma} (c-\#$branches of $\sigma$ in $T_i)$.

Take $\gamma_0=\omega^c$ and, for $i>0$:
\[
\gamma_i= \bigoplus_{\mathclap{\substack{ j<c \\ \emptyset\neq \sigma \in T_i \\  j = d_\sigma}}} \omega^j \cdot n_{\sigma, i}.
\]

One can check that: $\mathrm{MC} (\gamma_{i}) \leq E(x_i + c)$.

Notice that, by the absence of an $(\omega \cdot c)$-large subset of $X$: $\gamma_{i+1} < \gamma_i$ and $\gamma_{|X|} >0$.

This is a contradiction due to Lemma~\ref{lemma:411replacement}.
\qed
\begin{remark}
Any $\omega^{c+4}$-large $X>2$ is also $(\omega^{c+3}+\omega^3 + c+4)$-large.

\emph{Proof:} This follows from $\omega^{c+4}[x_0] \dots [x_{c+4}] >  (\omega^{c+3}+\omega^3 + c+4)$ with Lemmas~\ref{lemma:shiftlarge} and ~\ref{lemma:largerordinal},
\qed
\end{remark} 

\begin{lemma}\label{lemma:stepup}
Suppose $X$ is $(\omega^3 \cdot \omega^\alpha + \omega^3 + \max \{c ,\mathrm{MC}(\alpha)\}+3)$-large, then for every colouring $C\colon [X]^{d+1} \rightarrow c$ there exists $\alpha$-large ${\min}_d$-$C$-homogeneous subset of $X$. 
\end{lemma}
\emph{Proof:} Assume, for a contradiction, that the colouring $C\colon [X]^{d+1} \rightarrow c$ is such that it does not have $\alpha$-large ${\min}_d$-homogeneous subset of $X=\{ x_0 < \dots < x_{|X|-1}\}$. Define $T_0 \subset \dots \subset T_{|X|}$ as follows: $T_0= \{ \emptyset \}$ and 
\[
T_{i+1}= T_i \cup \{ \sigma \conc x_i\},
\]
where $\sigma \in T_i$ is the leftmost of maximum length such that $\{ \sigma_0 ,  \dots, \sigma_{\mathrm{lh}(\sigma )-1} ,x_i\}$ is ${\min}_d$-$C$-homogeneous. By construction, if $\sigma\conc y , \sigma\conc z \in T_i$ then there are $\sigma_{j_0} < \dots < \sigma_{j_{d-1}}$ with 
\[
C( \sigma_{j_0}  , \dots , \sigma_{j_{d-1}}, y ) \neq C(\sigma_{j_0}  , \dots , \sigma_{j_{d-1}}, z).
\]
So the number of branches of $\sigma$ is bound by the number of colourings $[\sigma_0, \dots , \sigma_{\lh (\sigma )-1} ]^d \rightarrow c$, which has upper bound $c^{2^{\sigma_{\mathrm{lh}(\sigma)-1}}}$. 

Define: $m_{\sigma, i}= c^{2^{\sigma_{\mathrm{lh}(\sigma)-1}}}-\#$branches of $\sigma$ in $T_i$.

By the comment directly above $m_{\sigma, i} \geq 0$.

Take $\gamma_0=\omega^\alpha$ and:
\[
\gamma_i= \bigoplus_{\emptyset \neq \sigma \in T_i} \omega^{\alpha[\sigma_0]\cdots [\sigma_{\mathrm{lh}(\sigma)-1}]}\cdot m_{\sigma, i}.
\]

One can check that: $\mathrm{MC} (\gamma_{i}) \leq E(x_i+ \max \{c, \mathrm{MC}(\alpha)\} +1)$.

We can see that, by the absence of ${\min}_d$-homogeneous $\alpha$-large subsets of $X$: $\gamma_i > \gamma_{i+1}$ and $\gamma_{|X|}>0$. 

This is a contradiction by Lemma~\ref{lemma:411replacement}.
\qed

Theorem~\ref{thm:ketonensolovay} can now be shown using induction on $d$ to prove:

If $X\geq  3$ is $\omega_{d}(\omega^{c+4}+d)$-large, then every colouring $C\colon [X]^{d+2} \rightarrow c$ has homogeneous $H\subseteq X$ of size $> \min H$.

Use Lemma~\ref{lemma:ketonensolovayd=0} for the base case and Lemma~\ref{lemma:stepup} in the induction step. Use Lemmas~\ref{lemma:shiftlarge} and ~\ref{lemma:largerordinal} to bridge the differences in largeness.

\begin{corollary}
Theorem~\ref{thm:ketonensolovay} is provable in $\RCA$. 
\end{corollary}

\begin{remark}
Above proof is also fine in $\RCAst$ if we fix a standard $d$.
\end{remark}

\section{A note on weakening the base theory}
We work in $\RCAst$ as defined in X.4 from \cite{sosoa} and elaborated in Section~2 from \cite{simpsonsmith}. Since every elementary function's existence is proven within $\RCAst$ we recommend Sections~2 and ~3 from Chapter~1 of \cite{rose} as background material. Defining our ordinals as in Section~\ref{section:ordinals} poses no problem, however note that,  in the proof in $\RCA$, we used implicitely that $(\alpha, x) \mapsto \alpha[x]$ is primitive recursive. 

As we are now working within the weaker system, we need this function to be elementary, which is non-obvious for nonstandard $d$. For example, if we encode the ordinals using prime numbers, as in \cite{simpsonsmith}, then we may need a nonstandard amount of iterations of the exponential function to determine the code of $\alpha[x]$, which is not available in $\RCAst$. To solve this problem, we will give an explicit encoding of ordinals which is consistent with the previous definitions and which will allow us to use bounded recursion to define fundamental sequences. 

Starting with the pairing map $j (x,y)= \frac{1}{2}(x+y+1)(x+y)+y$ and projections $j_1$ and $j_2$ for this map, we use:
\[
\pi^n_i (x)= \left\{ 
\begin{array}{ll}
j_1 j_2^{(n-i)}(x) 	& \textrm{if $1<i\leq n$}, \\
j_2^{(n)}(x)			& \textrm{if $i=1$}.
\end{array}
\right.
\]
Notice that:
\[
 j(\pi^n_n(x) , j(\pi^n_{n-1}(x), \dots , j(\pi^n_2(x), \pi^n_1(x)) \dots )))=x.
\]
Intuitively, the $\pi^n_i$ are the usual projection functions defined on alternatively coded $n$-tuples, for all (including non standard) $n$. 

\begin{definition}[Codes of ordinals]
Using bounded recursion, we define the codes of ordinals from $\mathcal{E}$ and relation $\prec$ on the codes as follows:
\begin{enumerate}
  \item $a$ is a code whenever $j_1(a)=n>0$ and $\pi^n_{1}(j_2(a)) \succeq \dots \succeq \pi^n_n (j_2(a))$ are codes.  
  \item Given codes $a,b$ with $n=j_1(a)$, $m=j_1 (b)$, $a'=j_2(a)$ and $b'=j_2(b)$, 
  
  $a \prec b$ if and only if:
  \begin{enumerate}
      \item $n<m$ and $\pi^n_i (a') = \pi^m_i (b')$ for all $0<i\leq m$, or
      \item there is $0<i \leq \min\{n , m\}$ with $\pi^{n}_j(a')=\pi^n_j(b')$ for all $0<j<i$ and $\pi^{n}_i(a')\prec \pi^n_i(b')$
  \end{enumerate}
\end{enumerate}
$0$ is the code for $0$, $w_0=j(1,0)$, $w_{i+1}=j(1, w_i)$ and $0 \preceq a$ for all codes of ordinals $a$. As we did with with the ordinals, we use $a_i=\pi^n_i(j_2(a))$.
\end{definition}

One can define ordinal addition, multiplication, natural (Hessenberg) sum, the Cantor Normal Form (coded version) and the maximum coefficient on the codes of the ordinals using bounded recursion. 

Using bounded recursion one can define a function $\mathrm{code}\colon \mathcal{E} \mapsto \N$ such that it preserves the order and operations on the ordinals. Furthermore $\mathrm{code}(\omega_d)=w_d$. 

Our next step is to define the fundamental sequences on the codes of ordinals. 

\begin{definition}
On the codes of ordinals, following the definition of fundamental sequences on ordinals, define: 
\[
a[x]= \left\{ 
\begin{array}{ll}
0 						& \textrm{if $a=0$}, \\
j(n-1, j_2^{(2)}(a)) 		& \textrm{if $n=j_1(a)>0$ and $a_n=0$}, \\
j(n+x-1, f(x, b, a) )			& \textrm{if $n=j_1(a)>0$, $m=j_1(a_n)>0$, $(a_n)_m=0$} \\
						& \textrm{and $b=j(m-1, j_2^{(2)}(a_n)) $},\\
j(n, j(a_n[x], j^{(2)}_2(a))	& \textrm{otherwise}, \\
\end{array}
\right.
\]
where $f(0, b, a)=j_2^{(2)}(a)$ and $f(i+1, b,a) = j(b, f(i, b,a))$.  
\end{definition}
\begin{lemma}
The functions $(a,x) \mapsto a[x]$ and $(a , \{x_0 < \dots < x_n\}) \mapsto a[x_0] \dots [x_n]$ are elementary.
\end{lemma}
\emph{Proof:} One can check:
\[
f(i, b,a) \leq 2^i (a+b+1)^{2i},
\]
hence, for $0<a\prec w_d$ and $x>0$:
\[
a[x] \leq 3^d (2^{x+1}a^{2x})^{2^d}.
\]
So:
\[
a[x_0] \dots [x_n] \leq (6a)^{d(2x_n+2)^{(n+1)(d+1)}}. 
\]
Therefore, these functions are elementary by bounded recursion. 
\qed
Taking care to use $\Delta^0_0$-induction every time induction was used, one can proceed with the proofs as described, using codes of ordinals instead of $\mathcal{E}$ where necessary. Simply observe that if a set is $\alpha$-large, then it is also $\mathrm{code}(\alpha)$-large, where $a$-largeness is similar to $\alpha$-largeness, but defined on the codes of ordinals instead of on the ordinals.

As an example of dealing with the induction steps, examine Lemma~\ref{lemma:fundbound} modified to the ordinal codes:
\begin{lemma}
If $w_d \succ a \succ b$ and $x>\MC(b)$, then $a[x] \succeq b$, where the inequality is strict if $a$ is a limit. 
\end{lemma}
\emph{Proof:} Given $x$, $a$ and $b$, use $\Delta^0_0$-induction on $d$ to prove the following:

If $w_d \succ a' \succ b'$, $a' \leq a$, $b'\leq b$, and $x>\MC(b)$, then $a'[x] \succeq b'$, where the inequality is strict if $a'$ is a limit.

To show that this can be expressed with a $\Delta^0_0$-formula, notice that the characteristic functions of $a \prec b$ and $a[x] \succeq b$ are elementary and use those functions as set parameter values.
\qed

\begin{corollary}
Theorem~\ref{thm:ketonensolovay} is provable in $\RCAst$.
\end{corollary}

\end{document}